\documentclass[12pt]{amsart}

\usepackage{amsfonts,amsmath,latexsym,amssymb,verbatim,amsbsy,amsthm}

\usepackage{graphicx}
\usepackage{grffile}  
\usepackage{nicefrac}
\usepackage{dsfont}
\usepackage{mathrsfs}
\usepackage{cancel}
\usepackage[normalem]{ulem}
\usepackage{color}
\usepackage{caption,subcaption}
\usepackage{url}
\usepackage[colorlinks=true, pdfstartview=FitV, linkcolor=blue,citecolor=red, urlcolor=blue]{hyperref}

\usepackage[dvipsnames]{xcolor}
\usepackage{relsize}

\renewcommand{\geq}{\geqslant}

\renewcommand{\leq}{\leqslant}

\newcommand{\be}{\begin{equation}}
\newcommand{\ee}{\end{equation}}

\numberwithin{equation}{section}
\numberwithin{figure}{section}
\numberwithin{table}{section}

\theoremstyle{plain}
\newtheorem{THEOREM}{Theorem}[section]

\newtheorem{theorem}[THEOREM]{Theorem}

\newtheorem{lemma}[THEOREM]{Lemma}

\theoremstyle{definition}

\theoremstyle{remark}
\theoremstyle{question}

\newtheorem{remark}[THEOREM]{Remark}

\newcommand{\bu}{{\mathbf u}}
\newcommand{\bw}{{\mathbf w}}

\newcommand{\R}{\mathbb{R}}

\newcommand{\lamb}{\lambda}
\newcommand{\Du}{\Delta \bar{u}}

\newcommand{\Dx}{\Delta x}
\newcommand{\Dy}{\Delta y}
\newcommand{\Dt}{\Delta t}
\newcommand{\Df}{\Delta f}
\newcommand{\Dbf}{\Delta {\mathbf f}}
\newcommand{\bp}{{\mathbf p}}
\newcommand{\Dg}{\Delta g}
\newcommand{\parp}{\partial_{x}^{+}}
\newcommand{\parm}{\partial_{x}^{-}}
\newcommand{\parx}{\partial_{x}}

\newcommand{\hf}{\frac{1}{2}}
\newcommand{\jph}{j+\hf}
\newcommand{\jmh}{j-\hf}

\newcommand{\nph}{n+\hf}
\newcommand{\baru}{\bar{u}}
\newcommand{\barbu}{\bar{\bu}}
\newcommand{\eps}{\epsilon}
\newcommand{\dx}{\textnormal{d}x}
\newcommand{\dt}{\textnormal{d}t}
\newcommand{\CFL}{{CFL} }
\newcommand{\vA}{\textnormal{vA}}

\DeclareMathOperator{\minmod}{minmod}

\title[Revisiting  high-resolution schemes with van Albada  slope limiter]{Revisting  high-resolution schemes\\with van Albada slope limiter}

\author{Jingcheng Lu}
\address{Department of Mathematics, University of Maryland, College Park}
\email{jlu1@umd.edu}

\author[Eitan Tadmor]{Eitan Tadmor*}
\address{Department of Mathematics and Institute for Physical Sciences \& Technology\newline \hspace*{0.45cm}University of Maryland, College Park}
\email{tadmor@umd.edu}
\date{\today}

\thanks{*Corresponding Author}
\thanks{\textbf{Acknowledgment.} Research was supported in part by ONR grant N00014-2112773.}



\dedicatory{\mbox{ }\newline{\large To R\'{e}mi Abgrall with friendship and appreciation}}

\begin{document}

\maketitle

\begin{abstract}
Slope limiters play an essential role in maintaining the non-oscillatory behavior of high-resolution methods for nonlinear conservation laws. The family of minmod limiters serves as the prototype example. Here, we revisit the question of non-oscillatory behavior of high-resolution central schemes in terms of the  slope limiter proposed by van Albada et. al. 1982. The  van Albada (vA) limiter is smoother near extrema, and consequently,  in many cases, it outperforms  the results obtained using the standard minmod limiter. In particular, we prove that the vA limiter ensures 1D TVD stability and  demonstrate that it yields noticeable improvement in computation of one- and two-dimensional systems.
\end{abstract}

\tableofcontents

\section{Introduction}
We revisit the class of high-resolution, non-oscillatory schemes for the approximate solution of nonlinear conservation laws. This class of schemes 
went through an intense period of development during the `80s and `90s.
By \emph{high-resolution} we refer to the class of schemes which are formally  second- or higher-order methods at all computational cells, except perhaps  for finitely many critical cells, thus breaking the first-order accuracy barrier of Godunov for monotone schemes, \cite{HHL76}. By \emph{non-oscillatory} we refer to  classes of schemes which satisfy certain type of stability bounds --- weaker than monotonicity yet strong enough to exclude spurious oscillations. In this context, we mention the canonical example of  Total Variation Diminishing (TVD) schemes of Harten \cite{Harten83,Sweby84,CO85,OT88,PD98}. Other classes of non-oscillatory stable schemes followed, including finite-volume, finite-element, Discontinuous Galerkin (DG) and spectral schemes; as examples we mention the Total-Variation  (TVB) bounded schemes \cite{Shu87,JSH90,Tadmor93,LSTZ08, WX19}, the One-Sided Lip Condition  (OSLC) bounded E-schemes  \cite{Osher84, Tadmor84}, and  entropy stable schemes \cite{Tadmor03, FMT12}. We also mention the notable class of Essentially Non-Oscillatory (ENO) and Weighted ENO schemes \cite{HEOC87,LOT94,Shu98,Shu20} and their sign property stability bound \cite{FMT13,FMT12,FR16}.\newline 
In this paper we focus on second-order Godunov-type methods. The key step of such  methods is the second-order accurate \emph{reconstruction}. To fix our notations, we consider one-dimensional,  piece-wise constant scalar numerical solution at  time level $t^n := n\Dt$, which is realized in terms  of its cell-averages, $\{\baru_\alpha^n= \baru(x_\alpha,t^n)\}$, over equi-spaced mesh $\{x_{\alpha} :  x_\alpha= \alpha{\Dx}\}$,
\begin{equation*}
    u^n = \sum^{}_{j} \baru^n_j {\mathds 1}_{[x_{\jmh}, x_{\jph}]}, \qquad  \baru(x_j,t^n) := \frac{1}{\Dx}\int^{x_{\jph}}_{x_{\jmh}}u(x,t^n) \dx, \qquad j\in {\mathbb Z}.
\end{equation*}
The key ingredient of Godunov-type schemes is a  piecewise linear  reconstruction, $u^{\Dx}(x,t^n) \approx u(x,t^n)$, which conserves the underlying cell averages, $\displaystyle \frac{1}{\Dx}\int^{x_{\jph}}_{x_{\jmh}}u^{\Dx}(x,t^n) \dx=\baru^n_j$. Thus, $u^{\Dx}(x,t^n)$ takes the form
\begin{equation}\label{eq:MUSCL interpolant}
\begin{split}
    u^{\Dx}(x,t^n) = \sum_{j} p^n_j(x) {\mathds 1}_{[x_{\jmh}, x_{\jph}]}, \qquad\
  p^n_j(x) = \baru^{n}_{j}+\frac{u'_j}{\Dx}(x-x_j),
\end{split}    
\end{equation}
 Here, $\displaystyle \frac{u'_j}{\Delta x}$ is an approximate slope that satisfies the second order of accuracy condition
\begin{equation}\label{eq:u slope}
    \frac{u'_j}{\Dx} = \frac{\partial}{\partial x}u(x_j,t^n)+{\mathcal O}(\Dx).
\end{equation}
The approximate slopes, $\{u'_j\}$, are reconstructed from the known cell averages $\{\baru^{n}_{j}\}$, so that they maintain second-order accuracy \eqref{eq:u slope}, while at the same time, they maintain a non-oscillatory property of the piecewise linear  reconstruction, $u^{\Dx}(x,t^n)$. A canonical  example for such reconstructed slopes is given by the so-called minmod limiter
\begin{subequations}\label{eqs:mionmod}
\begin{equation}
u'_j= \minmod(\Du^n_{\jmh}, \Du^n_{\jph}), \qquad \Delta u^n_{\jph} := u^n_{j+1}-u^n_j, 
\end{equation}
where $\minmod(\cdot,\cdot)$ is given by 
\begin{equation}\label{eq:minmod}
    \minmod(a,b) := \frac{sgn(a)+sgn(b)}{2}\cdot \min\{|a|,|b|\}.
\end{equation}
\end{subequations}
It follows that the corresponding reconstruction, $u^{\Dx}(x,t^n)$, satisfies the second-order accuracy at all by but the critical cells where $\frac{\partial}{\partial x}u(x_j,t^n)=0$, while satisfying the  TVD stability property
\begin{equation}\label{eq:TVD}
    ||u^{\Dx}(\cdot,t^n)||_{TV} \leq ||u^{n}||_{TV}:= \sum_j|\Delta u^n_{\jph}|.
\end{equation}
We now put this  second-order, TVD minmod-based reconstruction into action, in the context of approximate solution for one-dimensional nonlinear systems of conservation laws,
\begin{equation}\label{eq:CLS}
    \frac{\partial }{\partial t} \bu(x,t) +\frac{\partial }{\partial x} {\mathbf f}(\bu(x,t)) = 0, \qquad  \bu(x,t):\Omega\times \R_+ \mapsto \R^k, \quad \Omega\subseteq\R,
\end{equation}
 where  ${\mathbf f}$ is the $k$-vector of smooth flux functions. To this end, the reconstructed solution $\bu^{\Dx}(\cdot,t^n)$ is evolved to time $t^{n+1}=t^n+\Delta t$ by solving \eqref{eq:CLS} subject
 to the piecewise linear data, $\bu(\cdot,t^n)=\bu^{\Dx}(\cdot,t^n)$, and  the resulting solution, $\bu(x,t^{n+1})$ is then realized by its cell averages, $\bu^{n+1} = \sum^{}_{\alpha} \barbu^{n+1}_\alpha {\mathds 1}_{[x_{\alpha+\hf}, x_{\alpha-\hf}]}$. This completed the cycle of Godunov-type scheme: it consists of the  the reconstruction step  
 \addtocounter{equation}{1}
 \begin{equation}\tag*{(\theequation)${}_{\mathcal R}$}\label{eq:R}
 {\mathcal R}: \barbu^n_\alpha= \sum^{}_{\alpha} \barbu^n_\alpha {\mathds 1}_{[x_{\alpha-\frac{1}{2}}, x_{\alpha+\frac{1}{2}}]} \mapsto \bu^{\Dx}(\cdot,t^n);
\end{equation}
for simplicity, we assume that the reconstruction step is carried out by implementing the piece-wise linear reconstruction \eqref{eq:MUSCL interpolant},\eqref{eq:u slope}  component-wise\footnote{In general, one makes use of local characteristic decomposition.}. This then is followed by the evolution step, 
\begin{equation}\tag*{(\theequation)${}_{\mathcal E}$}\label{eq:E}
  {\mathcal E}(t^n,\Delta t): \bu^{\Dx}(\cdot,t^n) \mapsto \bu^{n+1}(\cdot,t^{n+1}),
 \end{equation}
   and ending with averaging 
\begin{equation}\tag*{(\theequation)${}_{\mathcal A}$}\label{eq:A}
   {\mathcal A}: \bu^{n+1}(\cdot,t^{n+1})\mapsto \bu^{n+1} = \sum^{}_{\alpha} \barbu^{n+1}_\alpha {\mathds 1}_{[x_{\alpha-\hf},x_{\alpha+\hf}]}.
\end{equation}
 It is well-known that shock discontinuities can --- and in most cases of generic data ---  will be formed at finite time with entropic solutions of \eqref{eq:CLS}.   Thus, one must accept that the underlying solution is discontinuous, and a high-resolution  numerical scheme should be designed so that it prevents the spurious (Gibbs) oscillations associated with such discontinuities, particularly when high-order  (order $>1$) accuracy is involved. Since averaging operator is inherently stable, it  is the primary role of a non-oscillatory  reconstruction  step, ${\mathcal R}$, to  keep the stability  compatible with the evolution step, e.g., \cite{Tadmor98}. Back to the scalar case, $k=1$, stability is interpreted in the sense of satisfying the TVD property; then,  \emph{if} the reconstruction is TVD, \eqref{eq:TVD} would imply the TVD stability of the overall scheme, 
 \[
 \begin{split}
\big\|\sum^{}_{\alpha} \baru^{n+1}_\alpha {\mathds 1}_{[x_{\alpha+\hf}, x_{\alpha-\hf}]}\big\|_{TV} &\leq \|{\mathcal E}(t^n,\Delta t)u^{\Dx}(\cdot,t^n)\|_{TV} \\
 & \leq \|u^{\Dx}(\cdot,t^n)\|_{TV} \\
  & \leq \big\|\sum^{}_{\alpha} \baru^n_j {\mathds 1}_{[x_{\alpha-\frac{1}{2}}, x_{\alpha+\frac{1}{2}}]}\big\|_{TV}  \leq \ldots \leq \big\|\sum^{}_{\alpha} \baru^0_\alpha {\mathds 1}_{[x_{\alpha-\frac{1}{2}}, x_{\alpha+\frac{1}{2}}]}\big\|_{TV}.
  \end{split}
 \]
 \begin{remark}[{\bf Staggered grids}]\label{rem:staggered} The last argument occupied  an arbitrary index $\alpha$ to emphasize that it may vary over the integers, $\alpha \leftrightarrow j\in {\mathbb Z}$, or over the half integers, $\alpha=j+\nicefrac{1}{2}$, or we can even allow staggered  grids with time-dependent shifts, $\alpha\mapsto\alpha^n=j+(\frac{1}{2})^{n-2[\nicefrac{n}{2}]}, \ j\in {\mathbb Z}$.
 \end{remark}
 The TV stability encoded in the successive  inequalities above, highlights the pivotal  role of the $\minmod$ limiter in securing the TVD stability of the overall scheme. 
 This line of argument extends to a general  class of non-oscillatory limiters 
 which take the form 
 \begin{equation*}
    u'_j = \psi(\Du^n_{\jmh}, \Du^n_{\jph}).
\end{equation*}
The corresponding Godunov-type scheme, $u^{n+1}(x)={\mathcal A}\,{\mathcal E}{\mathcal R}u^n(x)$ based on such $\psi$-limiter  amounts to a 5-cell stencil\footnote{To be precise, the resulting scheme is `essentially 3-point' stencil in the sense that the  corresponding minmod-based numerical flux  $F_{j+\frac{1}{2}}=F(u_{j-1},u_j,u_{j+1},u_{j+2})$ satisfies $F(\cdot,\bar{u},\bar{u},\cdot)=f(\bar{u})$, \cite{OT88}.}.  The essential feature of these limiters-based reconstructions is forming difference stencils in \emph{the direction of smoothness}, rather than in the direction of flow. Put differently, limiters include adaptive edge detectors   which are  necessary  in order to secure a non-oscillatory reconstruction, by avoiding `crossing' of the discontinuous data.
Higher order accuracy requires numerical derivatives with limiters which  occupy even wider stencils, e.g., \cite{CO85}
\[
\psi(...,  \Du^n_{j-\frac{3}{2}},\Du^n_{\jmh}, \Du^n_{\jph},\Du^n_{j+\frac{3}{2}}, ...),
\]
 which eventually led to the class of adaptive-based stencils in (W)ENO schemes, \cite{HEOC87,HO97,Shu98}. 
 Back to second-order resolution, the TVD criteria has profound influence on the development of high resolution numerical schemes, starting with the second order flux-limiter  MUSCL schemes \cite{vLeer79,vLeer84} and following the systematic framework for TVD stability offered in \cite{Harten83}. A variety of limiters have been proposed in the 1970s and were studied in the context of TVD limiters since 1980s. We mention two examples from the systematic study of second-order limiters  found in \cite{Sweby84}.  First is the class of $\minmod_\theta$ limiters
\[
u'_j= \minmod_\theta(\Du^n_{\jmh}, \Du^n_{\jph}), \qquad \Delta u^n_{\jph} := u^n_{j+1}-u^n_j, 
\]
where $\minmod_\theta$ is given by a one-parameter family of limiters
\begin{equation}\label{eq:theta minmod}
    \minmod_\theta(a,b) := \left\{\begin{array}{ll}
    \displaystyle s\cdot\min \Big(\theta |a|,\frac{|a+b|}{2},\theta |b|\Big), & \textnormal{if} \ sgn(a)=sgn(b):=s,\\
    0, & \textnormal{if} \ sgn(a)+sgn(b)=0,\end{array}\right. \qquad 1\leq \theta \leq 2.
\end{equation}
The case $\theta=1$ recovers the vanilla version of minmod \eqref{eq:minmod}.
As a second example we mention the van Albada (vA) limiter,
\begin{equation}\label{eq:VA eps limiter}
    \psi^{\eps}_{\vA}(a,b) = \frac{(a^{2}+\eps^{2})b+(b^{2}+\eps^{2})a}{a^{2}+b^{2}+2\eps^{2}}, \hspace{0.5em} \eps^2 = {\mathcal O}(\Dx^{3}),
\end{equation}
This limiter was first proposed by van Albada, van Leer and Roberts \cite{ALR82} in 1982 and has been successfully applied in computations (see e.g. \cite{MvL85,SA99,SPS05}), although a rigorous study of its TVD   stability was overshadowed by widespread applications of $\minmod$ limiter. A main aspect of this paper is to prove the TVD properties of vA limiter $\psi_{\vA}$. In fact, we demonstrate  that  when applied to the class of high resolution \emph{central schemes}, the vA limiter yields comparable  and in certain cases superior performance over the $\minmod$ limiter. 
This is attributed to the smoothness of the vA limiter at critical cells of local extrema\footnote{Specifically, the Lipschitz smoothness of $\psi_{\vA}(1,r)$ as a function of $r=\frac{\Du_{\jmh}}{\Du_{\jph}}$.},  in contrast to the $\minmod$ limiters, where $\minmod_\theta(\Du_{\jmh}, \Du_{\jph})=0$ whenever $\Du_{\jmh}\cdot \Du_{\jph}<0$.
Thus, the minmod-based polynomial reconstruction degenerates to first order piecewise constant interpolation at extrema cells, whereas the vA limiter still extracts a more accurate information in direction of smoothness, see \eqref{eq:VA-point},\eqref{eq:VA-cell} below.

\section{Brief Review of Second Order Central Schemes}
We distinguish between two classes of Godunov type schemes, depending 
on the averaging step in \ref{eq:A}. In one approach, averaging at $t=t^{n+1}$ is taken at the same grid as in $t=t^n$, with $\alpha \leftrightarrow j, \ j\in{\mathbb Z}$.
This leads to a the class of  \emph{upwind} schemes  which employ upstream-biased information to approximate spatial derivative.  The class of upwind scheme requires the solution of sequence of non-interacting Riemann solvers to realize the evolution step \ref{eq:E}. 
Wherever a second-order limiter is `turned off', the upwind scheme is reduced to the first-order Godunov scheme \cite{Godunov59} which evolves the piecewise constant approximate solution, based on cell averages $\{\barbu^n_j\}$, 
\begin{equation}\label{eq:godunov}
\begin{split}
    \barbu^{n+1}_{j} &= \barbu^{n}_{j}-\frac{1}{\Dx}\int \limits^{t^{n+1}}_{t^n} \big[{\mathbf f}(\bu^{\Dx}(x_{\jph},t))-{\mathbf f}(\bu^{\Dx}(x_{\jmh},t))\big] \textnormal{d}t = \barbu^{n}_{j}-\lambda\big( {\mathbf f}(\bu^{*}_{\jph})-{\mathbf f}(\bu^{*}_{\jmh})\big),
\end{split}
\end{equation}
where $\displaystyle \lambda := \frac{\Dt}{\Dx}$ is the fixed mesh ratio. The states $\bu^{*}_{j\pm\frac{1}{2}}$ are obtained with exact or approximate Riemann solvers, see e.g. \cite{EO81,HHL76,Roe81,Einfeldt88,TSS94}. When applied to systems of equations, both left- and right-running waves may exist at interfaces. Hence the local characteristic decomposition of the flux ${\mathbf f}$ is required for the upwind constructions, based on  the ``direction of the wind" identified by characteristic decomposition into local eigen-fields.

The class of \emph{central schemes}
keeps the averaging over staggered grid, $\alpha\leftrightarrow j+\nicefrac{1}{2}, \ j\in {\mathbb Z}$, so that it process the information from upstream and downstream in an averaged, ``central'' manner, \cite{NT90,Tadmor11}.
The evolution step in the class of central schemes requires simple quadrature; no Riemann solvers are needed. However, wherever a second-order limiter is `turned off', it is reduced to the diffusive Lax-Friedrichs  (LxF) scheme \cite{LxF54}:
\begin{equation}\label{eq:LxF}
    \barbu^{n+1}_{\jph} = \frac{\barbu^{n}_{j}+\barbu^{n}_{j+1}}{2}-\lambda({\mathbf f}(\bu^{n}_{j+1})-{\mathbf f}(\bu^{n}_{j})).
\end{equation}
 In contrast to the Godunov scheme \eqref{eq:godunov}, the LxF solver computes the cell averages on the staggered mesh, $\displaystyle \barbu(x_{\jph},t):= \frac{1}{\Dx}\int^{x_{j+1}}_{x_{j}}\bu(x,t)\dx$, recalling remark \ref{rem:staggered}. The integration of equation \eqref{eq:CLS} over the staggered control volume $[x_{j}, x_{j+1}]\times [t^n,t^{n+1}]$ yields
 \begin{equation}\label{eq:staggered cell average}
\barbu(x_{\jph},t^{n+1}) =\frac{1}{\Dx}\int_{x_{j}}^{x_{j+1}} \bu(x,t^{n}) \dx- \frac{1}{\Dx}\int^{t^{n+1}}_{t^n} [{\mathbf f}(\bu(x_{j+1},t))-{\mathbf f}(\bu(x_{j},t))] \dt.
 \end{equation}
It is clear that the LxF scheme is a first order approximation to \eqref{eq:staggered cell average}. Compared with the upwind framework, the advantage of central differencing is that the solution is smooth in the neighborhood of points $\{x_{j}\}$. Indeed, under a sufficiently small time step, $\Dt$, the Riemann waves initiating from  $x_{\jph}$ do not affect the adjacent cell centers within the time slab, $t\in[t^n, t^{n+1}]$. Hence the costly characteristic decompositions can be avoided. 

In spite of its simplicity, the main disadvantage of the LxF scheme is the excessive numerical dissipation which reduces its resolution at discontinuities and smooth extrema. The situation is dramatically improved by applying higher order piecewise polynomial interpolations, leading  to the second-order Nessyahu-Tadmor scheme\footnote{At the end, all fully-discrete second-order accurate schemes are one version of another of the Lax-Wendroff scheme.}. Still, wherever  a minmod limiter is `turned off',
the second-order central scheme is reduced to the difussive LxF scheme. 
 This is precisely why  the gain of a smooth vA limiter over the minmod limiter is more noticeable in the context of central schemes. 
We turn to briefly review the second order fully discrete and semi-discrete central schemes.

\subsection{The Nessyahu-Tadmor scheme}
The second-order fully discrete central scheme was first proposed in 1990 \cite{NT90}. The second-order of accuracy is obtained by combining the ${\mathcal R}\,{\mathcal E}{\mathcal A}$ steps. First, we evolve   a reconstructed piecewise linear solution $u(x,t^n)$  \eqref{eq:MUSCL interpolant}: integration of  the equation \eqref{eq:CLS} over the rectangle $[x_{j}, x_{j+1}]\times [t^n,t^{n+1}]$ yields
\begin{equation*}
    \barbu^{n+1}_{\jph} = \frac{1}{\Dx}\Big[\int^{x_{\jph}}_{x_{j}}\bp^n_{j}(x)\dx+\int^{x_{j+1}}_{x_{\jph}}\bp^n_{j+1}(x)\dx \Big]-\frac{1}{\Dx}\int^{t^{n+1}}_{t^n} \Big[{\mathbf f}(\bu^{\Dx}(x_{j+1},t))-{\mathbf f}(\bu^{\Dx}(x_{j},t))\Big] \dt.
    \end{equation*}
    Then, the time integral of the flux function can approximated by the midpoint rule at the expense of ${\mathcal O}(\Dt^3)$ local truncation error. This results in the Nessyahu-Tadmor (NT) scheme
\begin{equation}\label{eq:NT scheme}
    \barbu^{n+1}_{\jph} = \frac{\barbu^{n}_{j}+\barbu^{n}_{j+1}}{2}+\frac{1}{8}(\bu'_j-\bu'_{j+1})-\lambda[{\mathbf f}(\bu^{\nph}_{j+1})-{\mathbf f}(\bu^{\nph}_{j})],
\end{equation}
The midpoint values, $\bu^{\nph}_{j}$, are predicted by Taylor expansions
\begin{equation*}
    \bu^{\nph}_{j} = \barbu^n_{j}-\frac{\lambda}{2} {\mathbf f}'_j,
\end{equation*}
where the numerical derivative $f'_j$ satisfies the second order accuracy condition
\begin{equation}\label{eq:f slope}
    \frac{{\mathbf f}'_j}{\Dx} = \frac{\partial }{\partial x} {\mathbf f}(\barbu^n_j)+{\mathcal O}(\Dx).
\end{equation}
Different options of evaluating ${\mathbf f}'_j$ were proposed in \cite{NT90}. For example, the numerical derivative can be computed in terms of exact flux Jacobian, 
\begin{equation}\label{eq:f slope jacob}
    {\mathbf f}'_j = A(\barbu^n_{j})\bu'_j, \qquad A(\bu) := \frac{\partial {\mathbf f}(\bu)}{\partial \bu}.
\end{equation} 
Alternatively, one can apply Jacobian-free approximations when the characteristic decomposition is computationally expensive (even sometimes inaccessible), e.g.
\begin{equation}\label{eq:f slope jacob-free}
    {\mathbf f}'_j = \minmod(\Dbf_{\jmh}, \Dbf_{\jph}), \qquad \Dbf_{\jph} := {\mathbf f}(\barbu^n_{j+1})-{\mathbf f}(\barbu^n_j).
\end{equation}
Numerical results in \cite{JT98,LT98} have shown that the Jacobian-free version of the central scheme does not deteriorate the high resolution. When applied to systems of equations, the NT scheme inherits the simplicity of the LxF solver, i.e. the constructions can be extended in a component-wise manner without the use of Riemann solvers or characteristic decompositions. 

The multidimensional formulation of the NT scheme was obtained with similar integration procedures, details can be found in \cite{JT98,ASCM02,BTW04}. The order of accuracy can be further improved by applying higher order piecewise polynomial interpolants, see e.g. \cite{LT98,LPR99, LPR00,LPR02}. Related references can be found in \cite{CentPack}.

\subsection{The semi-discrete formulation}
It can be shown that the numerical dissipation of the NT scheme has order ${\mathcal O}((\Dx)^{4} /\Dt)$. In the convective problems where $\Dt\sim \Dx$, the NT scheme achieves higher resolution than the first order LxF scheme due to the reduced numerical viscosity. However, it is noticed that the NT scheme and its higher order extensions do not admit any semi-discrete limits, hence they are not appropriate for small time step computations or steady-state calculations. This motivated the development of semi-discrete central schemes. We sketch the derivations along the lines of \cite{KT00,KNP01}. 

We start with a piecewise polynomial approximation $\bu^{\Dx}(x,t^n)\approx u(x,t^n)$ of the form
\begin{equation*}
    \bu^{\Dx}(x,t^n) = \sum_{j} \bp^n_{j}(x){\mathds 1}_{[x_{\jmh}, x_{\jph}]}.
\end{equation*}
The polynomials $p^n_j(x)$ should have desired order of accuracy and conserve the cell averages $\barbu^n_j := \baru(x_j,t^n)$, i.e.
\begin{equation*}
    \frac{1}{\Dx} \int^{x_{\jph}}_{x_{\jmh}}\bp^n_j(x) \dx = \barbu^n_j.
\end{equation*}
We denote the reconstructed variables at $x_{\jph}$ from the left and the right by
\begin{equation}
    \bu^{-}_{\jph} := \bp^n_j(x_{\jph}), \qquad \bu^{+}_{\jph} := \bp^{n}_{j+1}(x_{\jph}).
\end{equation}
In particular, when the  MUSCL interpolant \eqref{eq:MUSCL interpolant} is applied we have 
\begin{equation}\label{eq:MUSCL reconstruction}
\bu^{-}_{\jph}: = \barbu^n_{j}+\frac{1}{2}\bu'_{j}, \hspace{1em} \bu^{+}_{\jph}: = \barbu^n_{j+1}-\frac{1}{2}\bu'_{j+1}.
\end{equation}

To obtain a semi-discrete formulation the key idea is to take cell averaging over narrower control volumes. Assume that the maximal forward and backward wave speeds at $x_{\jph}$ are estimated by $a^{+}_{\jph}\geq 0$ and $a^{-}_{\jph}\leq 0$. At the next time level $t^{n+1}:= t^n+\Dt$, the region influenced by the Riemann fan originating at $x_{\jph}$ is approximated with the interval, $I_{\jph} := [x^{-}_{\jph}, x^{+}_{\jph}]$, where 
\begin{equation*}
    x^{+}_{\jph}:= x_{\jph}+a^{+}_{\jph}\Dt, \hspace{1em}  x^{-}_{\jph}:= x_{\jph}+a^{-}_{\jph}\Dt.
\end{equation*}
 It is clear that under a sufficiently small time step, the non-smooth regions $I_{\jph}$ are separated due to the finite wave speeds. In this way, the width of Riemann fan is bounded by $(a^{+}_{\jph}-a^{-}_{\jph})\Dt$, in contrast to the fixed width $\Dx$ in the fully discrete counterpart.
 
 Now we reconstruct a non-oscillatory, conservative, piecewise polynomial interpolant at $t = t^{n+1}$,
\begin{equation}
    \widetilde{\bw}^{n+1}(x) = \sum_{j}\Big[\widetilde{\bw}^n_{j}(x){\mathds 1}_{[x^{+}_{\jmh}, x^{-}_{\jph}]}+\widetilde{\bw}^n_{\jph}(x){\mathds 1}_{ [x^{-}_{\jph}, x^{+}_{\jph}]}\Big].
\end{equation}
The polynomials $\widetilde{\bw}^{n+1}_{j}$ and $\widetilde{\bw}^{n+1}_{\jph}$ conserve the cell averages, $\bar{\bw}^{n+1}_{j}$ and $\bar{\bw}^{n+1}_{\jph}$, over the smooth and non-smooth domains. The values of $\bar{\bw}^{n+1}_{j}$ and $\bar{\bw}^{n+1}_{\jph}$ are computed by integrating \eqref{eq:CLS} over the rectangular domains, $[x^{+}_{\jmh}, x^{-}_{\jph}]\times[t^n,t^{n+1}]$ and $[x^{-}_{\jph}, x^{+}_{\jph}]\times[t^n,t^{n+1}]$, respectively.

Finally, the cell average at $t^{n+1}$ is computed with 
\[
\barbu^{n+1}_{j} = \frac{1}{\Dx}\int^{x_{\jph}}_{x_{\jmh}}\widetilde{\bw}^{n+1}(x)\dx.
\]
Passing the limit $\Dt\rightarrow0$ yields a semi-discrete conservative scheme
\begin{subequations}\label{eq:semi-discrete central scheme}
 \begin{equation}\label{eq:semi-discrete ode}
     \frac{\textnormal{d}}{\dt}\barbu_{j}(t) := \lim_{\Dt\rightarrow 0} \frac{\barbu^{n+1}_{j}-\barbu^n_{j}}{\Dt} = -\frac{{\mathbf F}_{\jph}(t)-{\mathbf F}_{\jmh}(t)}{\Dx},
 \end{equation}
where the numerical fluxes, ${\mathbf F}_{\jph}(t)$, are given by
\begin{equation}\label{eq:semi-discrete flux}
{\mathbf F}_{\jph}(t) := \frac{a^{+}_{\jph}{\mathbf f}(\bu^{-}_{\jph})-a^{-}_{\jph}{\mathbf f}(\bu^{+}_{\jph})}{a^{+}_{\jph}-a^{-}_{\jph}}+\frac{a^{+}_{\jph}a^{-}_{\jph}}{a^{+}_{\jph}-a^{-}_{\jph}}(\bu^{+}_{\jph}-\bu^{-}_{\jph}).
\end{equation}
\end{subequations}
There are different options to estimate the wave speeds $a^{\pm}_{\jph}$. The Kurganov-Tadmor (KT) scheme proposed in \cite{KT00}, uses the spectral radius of the Jacobian, $\displaystyle A(\bu):=\frac{\partial {\mathbf f}(\bu)}{\partial \bu}$,
 \begin{equation*}
     a^{+}_{\jph} = -a^{-}_{\jph} = a_{\jph} := \max\Big\{\rho(A(\bu^{-}_{\jph})), \rho(A(\bu^{+}_{\jph}))\Big\},
 \end{equation*}
where $\rho(\cdot)$ represents the spectral radius of a matrix. Then the numerical flux can be expressed in terms of the Rusanov flux,
\begin{equation}\label{eq:Rus flux}
    {\mathbf F}_{\jph}(t) = {\mathbf F}^{Rus}(\bu^{-}_{\jph}, \bu^{+}_{\jph}):= \frac{{\mathbf f}(\bu^{-}_{\jph})+{\mathbf f}(\bu^{+}_{\jph})}{2}-\frac{a_{\jph}}{2}(\bu^{+}_{\jph}-\bu^{-}_{\jph}).
\end{equation}
For strictly hyperbolic problems, the flux Jacobian $\displaystyle \frac{\partial {\mathbf f}}{\partial \bu}$ has $N$ distinct eigenvalues $\lambda_1<\cdots<\lambda_{N}$. The \emph{semi-discrete central-upwind scheme} version of KNP, \cite{KNP01}, employs a more accurate estimate of the wave speeds,
\begin{equation}\label{eq:hyper wave speeds}
\begin{split}
    & a^{+}_{\jph} = \max\Big\{\lambda_{N}(A(\bu^{-}_{\jph})), \lambda_{N}(A(\bu^{+}_{\jph})),0\Big\},\\
    & a^{-}_{\jph} = \min\Big\{\lambda_{1}(A(\bu^{-}_{\jph})), \lambda_{1}(A(\bu^{+}_{\jph})),0\Big\}.
\end{split}
\end{equation}
In this way, the scheme reduces the numerical dissipation by employing the HLL flux \cite{HLL83}.

In the computations, the set of ODEs \eqref{eq:semi-discrete central scheme} is integrated with an appropriate ODE solver. To preserve the overall high order of accuracy as well as the non-oscillatory property, one may apply a higher order \emph{Strong Stability Preserving} (SSP) Runge-Kutta method, see e.g. \cite{Shu88,GST01,SR02,GKS11}.

The semi-discrete central schemes retain the advantage of being Riemann-solver-free, hence the componentwise extension is allowed when solving systems of equations. In multidimensional problems, the second order constructions can be extended in a dimension-by-dimension manner. Third order extensions have also been derived with the help of piecewise parabolic interpolants, we refer the details to \cite{KT00,KL00}.

\section{Analysis of the Smooth van Albada Limiter}
In this section, we study the analytical properties of the smooth van Albada limiter \eqref{eq:VA eps limiter}. To get more insight, the limiter can be written as \cite{vLeer79} 
\begin{equation*}
    \psi^{\eps}_{\vA}(a,b) = \frac{a+b}{2}\Big(1-\frac{(a-b)^2}{a^2+b^2+2\eps^2}\Big).
\end{equation*}
In the smooth regions where $a\approx b$, the limiter tends to recover the second order central finite differencing, $\displaystyle \frac{a+b}{2}$. Across the discontinuities, however, the averaged slope is biased to the smallest value among the two one-sided slopes. These mechanisms are expected to ensure the second order accuracy and prevent the undesirable numerical oscillations. We will provide rigorous proofs for these advantages. For convenience of discussion, we will consider the simplified version without $\eps$
\begin{equation}\label{eq:VA limiter}
    \psi_{\vA}(a,b) = \frac{a^2 b+a b^2}{a^2+b^2}, \hspace{1em} a^2+b^2\neq 0.
\end{equation}
Indeed, as reported in \cite{ALR82}, the computed solutions are not sensitive to the specific value of $\eps$. 

\subsection{The second order accuracy}
We show that the van Albada limiter \eqref{eq:VA limiter} ensures a second order spatial accuracy. Back in the scalar framework, we first consider the reconstruction based on the \emph{point values}, $u_j := u(x_j)$. Denote the left- and the right- sided slopes by
\begin{equation*}
    \parm u_j:= \frac{u_{j}-u_{j-1}}{\Dx} \hspace{1em} \text{and} \hspace{1em} \parp u_j:= \frac{u_{j+1}-u_{j}}{\Dx}.
\end{equation*}
The numerical derivative constructed with the van Albada limiter \eqref{eq:VA limiter} is given by
\begin{equation*}
    \frac{u'_j}{\Dx} = \frac{(\parm u_{j})(\parp u_{j})^{2}+(\parp u_{j})(\parm u_{j})^{2}}{(\parp u_{j})^{2}+(\parm u_{j})^{2}}.
\end{equation*}
Taylor expansion yields
\begin{equation}\label{eq:VA-point}
    \frac{u'_j}{\Dx} = \parx u_{j} + \Big(\frac{1}{6}\parx^{3} u_{j}-\frac{1}{2}\frac{(\parx^{2} u_{j})^{2}}{\parx u_{j}}\Big)\Dx^{2}+{\mathcal O}(\Dx^{3}).
\end{equation}
The second order condition \eqref{eq:u slope} is satisfied in regions away from  critical points, $|\parx u_j|\gg \Dx$. In the framework of finite volume methods, the slope reconstruction is based on the \emph{cell averages}, $\displaystyle \baru_j := \frac{1}{\Dx}\int^{x_{\jph}}_{x_{\jmh}}u(x)\dx$. Then the Taylor expansion is slightly modified,
\begin{equation}\label{eq:VA-cell}
    \frac{u'_j}{\Dx} = \parx u_{j}+\Big(\frac{1}{6}\parx^{3}u_{j}+\frac{1}{24}\parx^{4}u_{j}-\frac{1}{2}\frac{\parx^{2}u_{j}}{\parx u_{j}}\Big)\Dx^{2}+{\mathcal O}(\Dx^{3}), \qquad |\parx u_j|\gg \Dx.
\end{equation}

For both point-value-based and cell-average-based reconstructions, the numerical derivative, $u'_j/\Dx$, approximates the exact derivative, $\parx u_j$, with order ${\mathcal O}(\Dx^2)$ in the non-critical regions. Near the critical points, the second order condition \eqref{eq:u slope} is satisfied up to $\parx u \approx \Dx$. Hence it is expected that the van Albada limiter introduces less dissipation than the minmod limiter at discontinuities and smooth extrema. Their numerical performances will be compared in section 4.

\subsection{The non-oscillatory property}
We are going to show that the van Albada limiter, as applied to the central schemes, generates non-oscillatory solutions. We begin with the following lemma. 

\begin{lemma}\label{lemma:VA estimate}
The approximate slope $u'_j = \psi_{\vA}(\Du_{\jmh},\Du_{\jph})$ satisfies the following estimates
\begin{subequations}\label{eq:VA estimate}
\begin{equation}\label{eq:VA bound}
\frac{1-\sqrt{2}}{2}\leq \frac{u_{j}'}{\Du_{j\pm\frac{1}{2}}} \leq \frac{1+\sqrt{2}}{2}, \hspace{1em} \end{equation}
\begin{equation}\label{eq:VA diff}
\frac{|u'_{j+1}-u'_{j}|}{|\Du_{\jph}|} \leq \sqrt{2}.
\end{equation}
\end{subequations}
\end{lemma}

\noindent \emph{Proof.} Denote $\displaystyle r := \frac{\Du_{\jmh}}{\Du_{\jph}}$, then we can write
\begin{equation*}
u'_{j} = \frac{r^{2}+r}{r^{2}+1} \Du_{\jph} = \frac{r+1}{r^{2}+1}\Du_{\jmh}.
\end{equation*}
 We will show that $\displaystyle \frac{1-\sqrt{2}}{2}\leq \frac{r^{2}+r}{r^{2}+1}, \hspace{0.2em}\frac{r+1}{r^{2}+1} \leq \frac{1+\sqrt{2}}{2}$. 
 
 Assume that $\displaystyle\frac{r^{2}+r}{r^{2}+1}$ is bounded from below and above by $m$ and $ M\in\R$ respectively,
\[
m\leq \frac{r^{2}+r}{r^{2}+1} \leq M, \hspace{0.5em} \forall r\in\R.
\]
This can be equivalently expressed as 
\begin{equation}
\left\{\begin{array}{l}
    (M-1)r^{2}-r+M \geq 0,\\
    \\
     (m-1)r^{2}-r+m \leq 0,
\end{array}\right. \qquad \forall r\in\R.
\end{equation}
The upper bound, $M$, should satisfy the conditions
\begin{equation*}
    \left\{\begin{array}{l}
         M>1,  \\
         \\
        \Delta =  1-4M(M-1)\leq 0. 
    \end{array}\right.
\end{equation*}
Solving these inequalities yields $\displaystyle M\geq \frac{1+\sqrt{2}}{2}$ and hence $\displaystyle \frac{1+\sqrt{2}}{2}$ gives the supremum of $\displaystyle \frac{r^2+r}{r^2+1}$. On the other hand, the lower bound, $m$, should satisfy 
\begin{equation*}
    \left\{\begin{array}{l}
         m<1,  \\
         \\
        \Delta =  1-4m(m-1)\leq 0. 
    \end{array}\right. 
\end{equation*}
These conditions yield $\displaystyle m\leq\frac{1-\sqrt{2}}{2}$ and hence $\displaystyle \frac{1-\sqrt{2}}{2}$ gives the infimum of $\displaystyle \frac{r^2+r}{r^2+1}$. Combining these results, we have
\begin{equation*}
    \frac{1-\sqrt{2}}{2}\leq \frac{u'_j}{\Du_{\jph}} = \frac{r^2+r}{r^2+1} \leq \frac{1+\sqrt{2}}{2}.
\end{equation*}
With similar argument, we can show that $\displaystyle \frac{1-\sqrt{2}}{2}\leq\frac{r+1}{r^{2}+1} \leq \frac{1+\sqrt{2}}{2}$ and therefore
\begin{equation*}
    \frac{1-\sqrt{2}}{2}\leq \frac{u'_j}{\Du_{\jmh}} = \frac{r+1}{r^2+1} \leq \frac{1+\sqrt{2}}{2}.
\end{equation*}
This completes the proof of \eqref{eq:VA bound}. The second estimate \eqref{eq:VA diff} is a direct consequence of \eqref{eq:VA bound}.
$\hfill \square$

With the help of Lemma \ref{lemma:VA estimate}, we discuss the non-oscillatory property of the fully discrete and the semi-discrete schemes separately. 

\subsubsection{The fully discrete scheme}
We consider the fully discrete NT scheme \eqref{eq:NT scheme}. The scheme can be written in a conservative form \cite{NT90}\footnote{Recall that $\lambda$ is the fixed mesh-ratio, $\lambda=\frac{\Delta t}{\Delta x}$, which should be distinguished from the eigenvalues in \eqref{eq:hyper wave speeds}, $\lambda_1<\ldots <\lambda_N$.}
\begin{equation}
\label{eq:NT conserv}
\baru_{\jph}^{n+1} = \frac{\baru^{n}_{j}+\baru^{n}_{j+1}}{2}-\lamb(g_{j+1}-g_{j}), \qquad \lambda=\frac{\Delta t}{\Delta x},
\end{equation}
with the so-called modified numerical flux, $g_j$,  given by
\begin{equation}\label{eq:NT modified flux}
g_{j} = f(u^{\nph}_{j})+\frac{1}{8\lamb}u'_{j}.
\end{equation}
The NT scheme combined with the van Albada limiter \eqref{eq:VA limiter} can be shown to be TVD following similar approach in \cite{OT88,NT90}. The proof is based on the following lemma.

\begin{lemma}\label{lemma:NT TVD condition} (\cite{NT90})
The scheme \eqref{eq:NT scheme} is TVD if the numerical flux, $g_j$, satisfies the following generalized \CFL condition
\begin{equation}\label{eq:NT TVD condition}
    \lambda \Big|\frac{\Dg_{\jph}}{\Du^n_{\jph}}\Big|\leq\frac{1}{2}, \hspace{1em} \Dg_{\jph} := g_{j+1}-g_{j}.
\end{equation}
\end{lemma}
\noindent \emph{Proof.} By \eqref{eq:NT conserv}, the difference $\baru^{n+1}_{\jph}-\baru^{n+1}_{\jmh}$ amounts to
\begin{equation*}
    \baru^{n+1}_{\jph}-\baru^{n+1}_{\jmh} = \Big(\frac{1}{2}-\lamb\frac{\Dg_{\jph}}{\Du_{\jph}}\Big)\Du_{\jph}+    \Big(\frac{1}{2}+\lamb\frac{\Dg_{\jmh}}{\Du_{\jmh}}\Big)\Du_{\jmh}.
\end{equation*}
Condition \eqref{eq:NT CFL} implies that the coefficients in the parenthesis are positive. Hence,
\begin{equation*}
\begin{split}
 \sum_{j}|\baru^{n+1}_{\jph}-\baru^{n+1}_{\jmh}| & \leq \sum_{j} \Big(\frac{1}{2}-\lamb\frac{\Dg_{\jph}}{\Du_{\jph}}\Big)|\Du_{\jph}|+    \Big(\frac{1}{2}+\lamb\frac{\Dg_{\jmh}}{\Du_{\jmh}}\Big)|\Du_{\jmh}|\\
& = \sum_{j} \Big(\frac{1}{2}-\lamb\frac{\Dg_{\jph}}{\Du_{\jph}}\Big)|\Du_{\jph}|+    \Big(\frac{1}{2}+\lamb\frac{\Dg_{\jph}}{\Du_{\jmh}}\Big)|\Du_{\jph}|\\
& = \sum_{j}|\Du_{\jph}|.
\end{split}
\end{equation*}
We conclude $TV(\baru^{n+1})\leq TV(\baru^n)$. $\square$

\bigskip\noindent
Equipped with Lemma \ref{lemma:VA estimate} and \ref{lemma:NT TVD condition}, we prove the NT scheme is TVD.

\begin{theorem}[{\bf TVD stability of NT scheme with vA limiter}]\label{thm:NT TVD}
Consider the scalar NT scheme  \eqref{eq:NT conserv} using numerical slope $u'_j$  constructed with the van Albada limiter \eqref{eq:VA limiter}, and let the flux numerical derivative be chosen by $f'_j = a(\baru^n_j)u'_j$ with $a(u) := f'(u)$. Assume that the following \CFL condition is satisfied
\begin{equation}\label{eq:NT CFL}
    \lambda\max_{u}|a(u)|\leq \frac{-1+\sqrt{\frac{3}{2}+\frac{3\sqrt{2}}{4}}}{1+\sqrt{2}} \approx 0.24.
\end{equation}
Then the NT scheme \eqref{eq:NT conserv} is TVD.
\end{theorem}

\noindent\emph{Proof.} Denote $\displaystyle \CFL:= \lambda\max_{u}|a(u)|$. By \eqref{eq:NT modified flux} we have
\begin{equation}\label{eq:NT decomp}
\lamb\Big|\frac{\Dg_{\jph}}{\Du_{\jph}}\Big|\leq \underbrace{\lamb\Big|\frac{f(u^{\nph}_{j+1})-f(u^{\nph}_{j})}{u^{\nph}_{j+1}-u^{\nph}_{j}}\Big|}_{I} 
\underbrace{\Big|\frac{u^{\nph}_{j+1}-u^{\nph}_{j}}{\Du_{\jph}}\Big|}_{II}+
\underbrace{\frac{1}{8}\Big|\frac{\Delta u'_{\jph}}{\Du_{\jph}}\Big|}_{III}.
\end{equation}

By definition, the first term on the right of \eqref{eq:NT decomp} is bounded by
\begin{equation}\label{eq:NT est I} 
    \lamb\left|\frac{f(u^{\nph}_{j+1})-f(u^{\nph}_{j})}{u^{\nph}_{j+1}-u^{\nph}_{j}}\right| \leq \lambda\max_u |a(u)| := \CFL.
\end{equation}
Notice that \eqref{eq:VA bound} implies $\displaystyle \Big|\frac{u'_j}{\Du_{j\pm\frac{1}{2}}}\Big|\leq \frac{1+\sqrt{2}}{2}$, we can estimate the second term on the right of \eqref{eq:NT decomp},

\begin{equation}\label{eq:NT est II}
\begin{split}
\Big|\frac{u^{\nph}_{j+1}-u^{\nph}_{j}}{\Du_{\jph}}\Big|& \leq 1+\frac{\lamb}{2}\Big|\frac{\Df'_{\jph}}{\Du_{\jph}}\Big|
 \leq 1+\frac{\lamb}{2} \frac{|a_{j+1}u_{j+1}'|+|a_{j}u'_{j}|}{|\Du_{\jph}|}\\
& \leq 1+\frac{\lamb\max_{u}|a(u)|}{2} \cdot \frac{|u_{j+1}'|+|u'_{j}|}{|\Du_{\jph}|}  \leq 1+\frac{1+\sqrt{2}}{2}\CFL .
\end{split}    
\end{equation}
As a direct consequence of \eqref{eq:VA diff}, the third term on the right of \eqref{eq:NT decomp} does not exceed
\begin{equation}\label{eq:NT est III}
    \frac{1}{8}\Big|\frac{\Delta u'_{\jph}}{\Du_{\jph}}\Big|\leq\frac{\sqrt{2}}{8}.
\end{equation}
Using \eqref{eq:NT est I}, \eqref{eq:NT est II} and \eqref{eq:NT est III}, the TVD condition \eqref{eq:NT TVD condition} boils down to 
\begin{equation*}
    \CFL\Big(1+\frac{1+\sqrt{2}}{2}\CFL\Big)+\frac{\sqrt{2}}{8}\leq \frac{1}{2},
\end{equation*}
which in turn  recovers the CFL condition \eqref{eq:NT CFL}. $\square$

\begin{remark}
The TVD stability is not sharp in the sense that the \CFL condition  \eqref{eq:NT CFL} is to restrictive, and serves here only as a theoretically sufficient bound. The numerical experiments reported in section \ref{sec:numerics}, indicate that the computed NT solutions with vA limiter remain non-oscillatory as long as $\displaystyle \CFL\leq \nicefrac{1}{2}$.
\end{remark}

\subsubsection{The semi-discrete scheme}
The TVD stability of the van Albada limiter as applied to the non-staggered semi-discrete schemes \eqref{eq:semi-discrete central scheme} can be shown following the lines of \cite{Osher85,Tadmor88}. The key ingredient is to represent the scheme in an appropriate incremental form which meets a certain positive condition. We prove the following result with the help of Lemma \ref{lemma:VA estimate}.

\begin{theorem}[{\bf TVD stability of semi-discrete MUSCL scheme with vA limiter}]\label{thm:MUSCL TVD}
Let the states $u^{\pm}_{\jph}$ be computed from the MUSCL reconstruction \eqref{eq:MUSCL reconstruction}, with the numerical slope $u'_j$ constructed with the van Albada limiter \eqref{eq:VA limiter}. Consider a generalized scalar MUSCL scheme
\begin{equation}\label{eq:generalized MUSCL scheme}
    \frac{d}{dt}\baru_j(t) = -\frac{1}{\Dx}[F(u^{-}_{\jph}, u^{+}_{\jph})-F(u^{-}_{\jmh}, u^{+}_{\jmh})], 
\end{equation}
and let the numerical flux $F(\cdot,\cdot)$ be monotone and consistent,
\begin{subequations}
\begin{equation}\label{eq:flux monotone}
    \frac{\partial}{\partial u} F(u,v) \geq 0 \quad \text{and} \quad \frac{\partial}{\partial v} F(u,v) \leq 0, 
\end{equation}
\begin{equation}\label{eq:flux consistency}
    F(u,u) = f(u).
\end{equation}
Then the scheme \eqref{eq:generalized MUSCL scheme} is TVD,
\begin{equation*}
    \frac{d}{dt} [TV(\baru(t))] \leq 0.
\end{equation*}
\end{subequations}

\end{theorem}

\noindent\emph{Proof.} We argue along the lines of \cite{Tadmor88}, expressing the scheme \eqref{eq:generalized MUSCL scheme} in the incremental form
\begin{equation}\label{eq:semi-discrete incremental}
\frac{d}{dt}\baru_{j} = -\frac{1}{\Dx_{j}}C_{\jmh}\Du_{\jmh}+\frac{1}{\Dx_{j}}D_{\jph}\Du_{\jph},
\end{equation}
where 
\begin{equation}\label{eq:incremental coeffs}
\begin{split}
&C_{\jmh} = \frac{1}{\Du_{\jmh}}[F(u^{-}_{\jph},u^{+}_{\jmh})-F(u^{-}_{\jmh}, u^{+}_{\jmh})], \\ 
&D_{\jph} = -\frac{1}{\Du_{\jph}}[F(u^{-}_{\jph}, u^{+}_{\jph})-F(u^{-}_{\jph},u^{+}_{\jmh})].
\end{split}
\end{equation}
We denote
\begin{equation*}
\begin{split}
&s_{\jph} := sign(\Du_{\jph}) = \left\{\begin{array}{cc}
1     &  \Du_{\jph}>0\\
\pm 1    & \Du_{\jph}=0\\
-1     & \Du_{\jph}<0,
\end{array}\right.  \\
&\chi_{j} := 1-s_{\jmh}s_{\jph} = \left\{\begin{array}{ll}
2,     &  \Du_{\jmh} \cdot \Du_{\jph}\leq 0\\
0,     & \Du_{\jmh} \cdot \Du_{\jph}> 0.
\end{array}\right.
\end{split}
\end{equation*}

Forward differencing of \eqref{eq:semi-discrete incremental} gives
\begin{equation}\label{eq:incremental forward diff}
\begin{split}
    \frac{d}{dt}\Du_{\jph} &= \Big(\frac{1}{\Dx_{j+1}}D_{j+\frac{3}{2}}\Du_{j+\frac{3}{2}}-\frac{1}{\Dx_{j+1}}C_{\jph}\Du_{\jph}\Big)\\
    &\qquad -\Big(\frac{1}{\Dx_{j}}D_{\jph}\Du_{\jph}-\frac{1}{\Dx_{j}}C_{\jmh}\Du_{\jmh}\Big).
\end{split}
\end{equation}
Multiply \eqref{eq:incremental forward diff} by $s_{\jph}$ and sum by parts, we have
\begin{equation*}
\begin{split}
\frac{d}{dt}TV(\baru(t)) & = \sum_{j} s_{\jph} \frac{d}{dt}\Du_{\jph}\\
& = -\sum_{j} \frac{1}{\Dx}[(s_{\jmh}-s_{\jph})C_{\jmh}\Du_{\jmh}+(s_{\jph}-s_{\jmh})D_{\jph}\Du_{\jph}]\\
& = -\sum_{j} \frac{\chi_{j}}{\Dx}[C_{\jmh}|\Du_{\jmh}|+D_{\jph}|\Du_{\jph}|].
\end{split}    
\end{equation*}
Here we use the property that $s^{2}_{\jph} \equiv 1$. Since $\chi_j\geq 0$, it suffices to show that $C_{\jmh}$, $ D_{\jph}\geq 0$.
Indeed, we find that 
\begin{equation}
\begin{split}
&C_{\jmh} = \partial_{u}F(\tilde{u}_1, u^{+}_{\jmh})\Big(1+\frac{u'_j-u'_{j-1}}{2\Du_{\jmh}}\Big), \\
&D_{\jph} = -\partial_{v}F(u^{-}_{\jph}, \tilde{u}_2)\Big(1-\frac{u'_{j+1}-u'_{j}}{2\Du_{\jph}}\Big),
\end{split}
\end{equation}
where $\tilde{u}_1$ is between $u^{-}_{\jmh}$ and $u^{-}_{\jph}$, $\tilde{u}_2$ is between $u^{+}_{\jmh}$ and $u^{+}_{\jph}$. From \eqref{eq:VA diff} we know that the terms in the parentheses are positive. Hence the non-negativity of $C_{\jmh}$ and $D_{\jph}$ follows from \eqref{eq:flux monotone}.
$\square$

\medskip
The semi-discrete central scheme \eqref{eq:semi-discrete central scheme} can be considered a generalized MUSCL scheme, the associated numerical flux \eqref{eq:semi-discrete flux} can expressed in the form $F_{\jph}(t) := F(u^{-}_{\jph}, u^{+}_{\jph})$, with
\begin{equation*}
    F(u,v) = \frac{a^{+}f(u)-a^{-}f(v)}{a^{+}-a^{-}}+\frac{a^{+}a^{-}}{a^{+}-a^{-}}(v-u).
\end{equation*}
The monotonicity of $F(\cdot,\cdot)$ is clear since $a^- \leq f'(u), f'(v) \leq a^+$. Hence the TVD stability of the van Albada limiter as applied to the semi-discrete central scheme \eqref{eq:semi-discrete central scheme} follows from Theorem \ref{thm:MUSCL TVD}.

\begin{remark}
The semi-discrecte scheme \eqref{eq:semi-discrete central scheme} (or equivalently \eqref{eq:generalized MUSCL scheme}) advanced with foward Euler time-stepping reads
\begin{equation}\label{eq:MUSCL FE}
    \baru^{n+1}_j = \baru^n_j-\lambda(F(u^{-}_{\jph}, u^{+}_{\jph})-F(u^{-}_{\jmh}, u^{+}_{\jmh})),
\end{equation}
where $u^{+}_{\jph} := \baru^{n}_j-\frac{1}{2}u'_{j+1}$ and $u^{-}_{\jph} := \baru^{n}_j+\frac{1}{2}u'_{j}$. With similar incremental representation in Theorem \ref{thm:MUSCL TVD}, it can be proved that the van Albada limiter as applied to the fully discrete non-staggered scheme \eqref{eq:MUSCL FE} generates TVD solutions. Applying a higher order SSP time integration method would achieve a higher order of accuracy in time without deteriorating the stability.
\end{remark}

\section{Numerical Experiments}\label{sec:numerics}
In this section, we examine the performance of vA limiter \eqref{eq:VA eps limiter} with the small\footnote{Even smaller than the one indicated in \eqref{eq:VA eps limiter}; this does not seem to affect the results.} bias, $\eps = (\Dx)^3$,  compared vs. the minmod limiter \eqref{eq:theta minmod} which  is taken, unless otherwise stated, in its vanilla version $\theta=1$. The advantage of the vA limiter over minmod${}_1$ is apparent  for both --- the fully discrete NT scheme \eqref{eq:NT scheme} and the semi-discrete central-upwind  version of KNP, \eqref{eq:semi-discrete central scheme}, \eqref{eq:hyper wave speeds}.

In all the simulations of the NT scheme, the numerical slope of flux, ${\mathbf f}'_j$ is evaluated with the exact flux Jacobian \eqref{eq:f slope jacob}, and the CFL number is taken to be 0.45. The implementation of the semi-discrete scheme requires a high order time discretization. We use the third order, explicit SSPRK3 method for the time integration of method-of-lines ODEs \eqref{eq:semi-discrete ode} with the CFL number chosen to be 0.7. To better ensure the robustness in the strong shock problems, the MUSCL reconstruction for systems of equations is performed on the characteristic variables, e.g., \cite[(4.7)-(4.11)]{NT90}, \cite[Procedure 2.8]{Shu98}\footnote{Here we use standard characteristic decomposition, $\Delta\bu_{j+\hf}=\sum_k \widehat{\alpha}^k_{j+\hf}\widehat{{\mathbf R}}^k_{j+\hf}$ where 
$\{\widehat{{\mathbf R}}^k_{j+\hf}\}_k$ is the  eigensystem of an intermediate Jacobian: for convenience,  $\widehat{{\mathbf R}}^k_{j+\hf}$ are computed as the eigensystem at the  arithmetic average, $A\big(\hf(\bu_j+\bu_{j+1})\big)$, rather than the Roe average,  $\Delta {\mathbf f}_{j+\hf}=A^{\textnormal{Roe}}_{j+\hf}\Delta \bu_{j+\hf}$, as it does not seem to affect the MUSCL results in this context.}

\subsection{One-dimensional linear advection equation}\label{sec:advection}
\begin{figure}[h!]
    \centering
    \includegraphics[scale = 0.4]{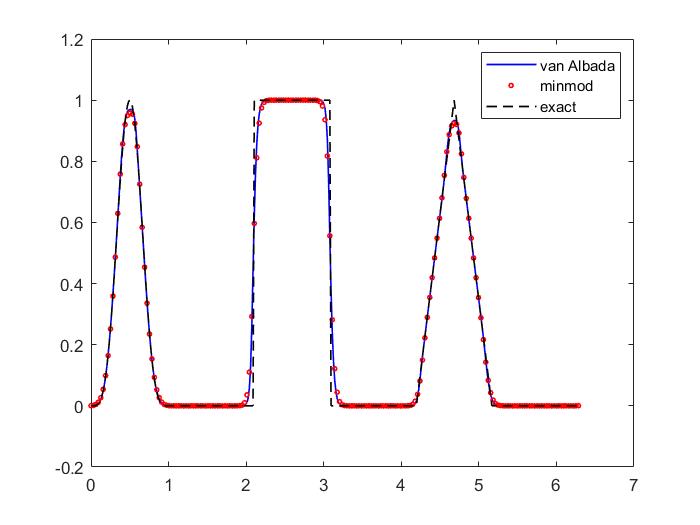}
    \caption{Linear advection. NT scheme, $N = 400$, periodic BCs, $t = 2\pi$.}
    \label{fig:multiwaves NT}
\end{figure}

\begin{figure}[h!]
    \centering
    \includegraphics[scale = 0.4]{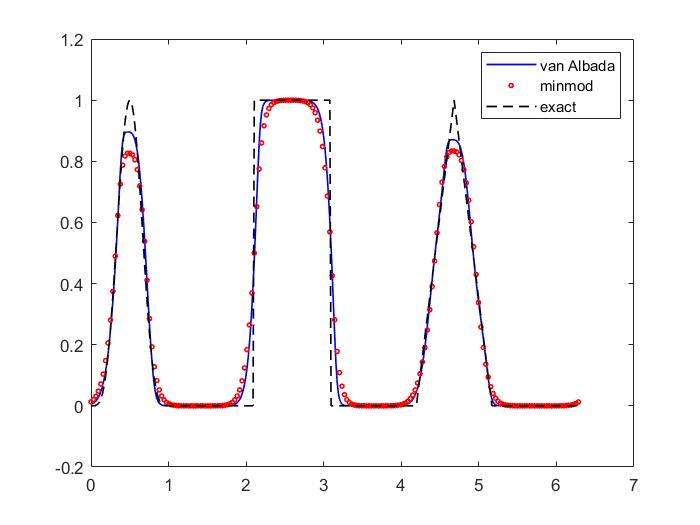} 
    \caption{Linear advection. Semi-discrete central-upwind scheme, $N = 400$, periodic BCs, $t = 2\pi$.}
    \label{fig:multiwaves csup}
\end{figure}

We consider the one-dimensional linear advection equation
\begin{equation*}
    \frac{\partial u}{\partial t}+\frac{\partial u}{\partial x} = 0, \qquad 0\leq x\leq 2\pi,
\end{equation*}
subject to initial data,
\[
u(x,0)\left\{\begin{array}{ll} 
 = \sin^4(\pi x) & 0\leq x \leq 1, \\
  \equiv 1 & 2.09 \leq x \leq 3.09, \\
 = \textnormal{hat function with height} = 1 & 4.18\leq x \leq 5.18,\\
 \equiv 0 & \textnormal{elsewhere in} \ (1,2\pi),
\end{array}\right.
\]
and  $2\pi$ periodic boundary conditions (BCs). We compute the solution after one-period of revolution using  $N = 400$ equi-spaced mesh grids. Figure \ref{fig:multiwaves NT} displays the solutions computed with the NT scheme. Both the vA limiter and the minmod limiter generate non-oscillatory solutions with comparable resolution at the shocks and smooth extrema. Figure \ref{fig:multiwaves csup} compares the solutions computed with different limiters as applied to the semi-discrete central-upwind scheme. In this case, applying the smoother vA limiter  improves the noticeable `clipping' of the minmod${}_1$ that occurred at the peaks of the sinusoidal wave and the triangular wave.

\begin{figure}[h!]
    \centering
    \includegraphics[scale = 0.4]{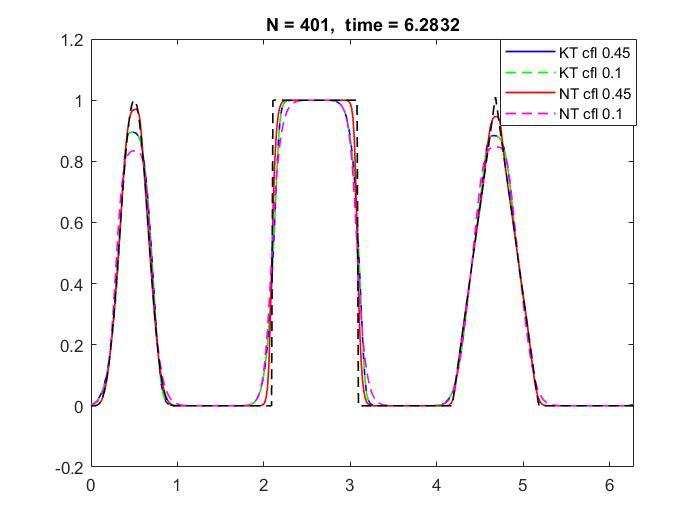} 
    \caption{Linear advection. NT vs. semi-discrete KT scheme with vA limiter with $N = 400$ mesh points,  periodic BCs at $t = 2\pi$.}
    \label{fig:multiwaves csupCFL}
\end{figure}

\begin{remark}[{\bf On the superior resolution of NT with vA limiter}]\label{rem:NTvA} It is instructive to compare the NT results in Figure \ref{fig:multiwaves NT} which are found to be superior to the corresponding results of the central-upwind KT scheme in Figure \ref{fig:multiwaves csup}. In particular, a direct comparison in Figure \ref{fig:multiwaves csupCFL} shows that the `clipping'
phenomenon in the three peaks is clearly stronger in the KT computation, whereas the NT scheme produces sharper peaks. We recall that the NT with minmod limiter is reduced at peaks to a diffusive Lax-Friedrichs scheme. Thus, the vA limiter helps  improving the NT resolution at peaks.   Of course, the truncation error of NT schemes of order ${\mathcal O}((\Delta x)^2/\lambda)$ implies that the NT resolution decreases for small(-er) CFL. But for near optimal CFL=0.45, the resolution of NT  scheme with vA limiter offers a competitive alternative to central-upwind schemes. This conclusion will be confirmed with \emph{linear} contact waves studied in section \ref{sec:1DEuler} below.
\end{remark}

\vspace{0.5in}

\newpage
\subsection{One-dimensional Burgers' equation}
\begin{figure}[h!]
     \centering
     \includegraphics[scale = 0.4]{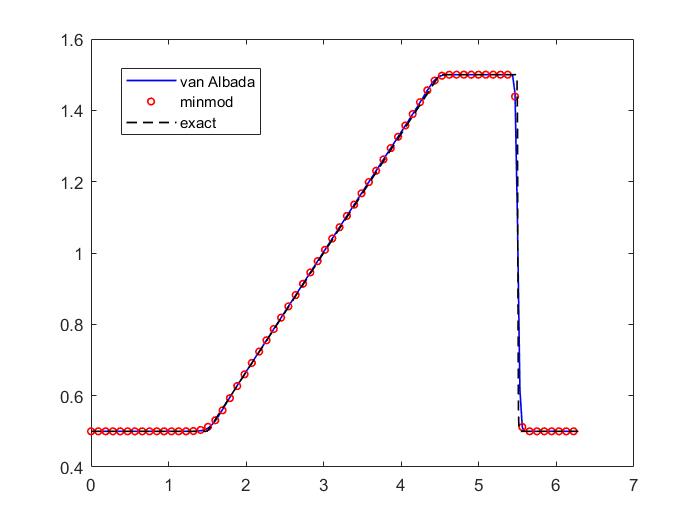}
     \caption{1D Burger's equation, NT scheme, $N = 200$, $t = 2$.}
     \label{fig:exp comp NT}
 \end{figure}

\begin{figure}[h!]
     \centering
     \includegraphics[scale = 0.4]{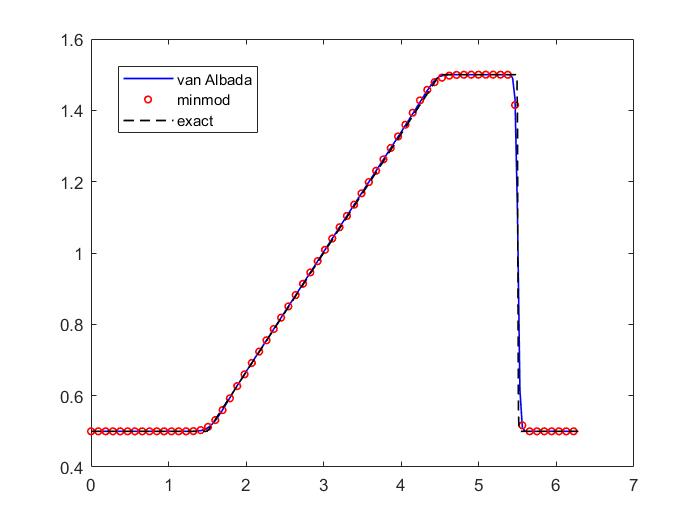}
     \caption{1D Burgers' equation, semi-discrete central-upwind scheme, $N = 200$, $t = 2$.}
     \label{fig:exp comp csup}
 \end{figure}

We show the results for the one-dimensional Burgers' equation
\begin{equation*}
    \frac{\partial u}{\partial t}+\frac{\partial}{\partial x}\Big(\frac{u^2}{2}\Big) = 0
\end{equation*}
 over the domain $[0, 2\pi]$. In this test case, the solution comprises of an expansion wave and a compression wave. Figures \ref{fig:exp comp NT},  \ref{fig:exp comp csup} present the solutions computed with $N = 200$ mesh grids at $t = 2$. In both fully discrete and semi-discrete test cases, the solutions generated with the van Albada limiter capture the expansion wave and the shock wave with sharp resolution and without generating numerical oscillations. The results based on vA limiter are essentially comparable to those based on the minmod limiter.

\newpage
 \subsection{One-dimensional Euler equations}\label{sec:1DEuler}
 We consider the one-dimensional Euler equations of gas dynamics
 \begin{equation*}
     \frac{\partial}{\partial t} \left[\begin{array}{c}
          \rho  \\
          \rho u\\
          E
     \end{array}\right]
     +\frac{\partial}{\partial x}\left[\begin{array}{c}
          \rho u  \\
           \rho u^2+p\\
           (E+p)u
     \end{array}\right] = 0,
 \end{equation*}
where $\rho$, $u$, $p$ are the density, velocity and pressure, $E$ is the total energy per unit volume given by
\[
E = \frac{p}{\gamma-1}+\frac{\rho u^2}{2}, \quad \gamma = 1.4.
\]

\begin{figure}[h!]
    \centering
    \includegraphics[scale = 0.4]{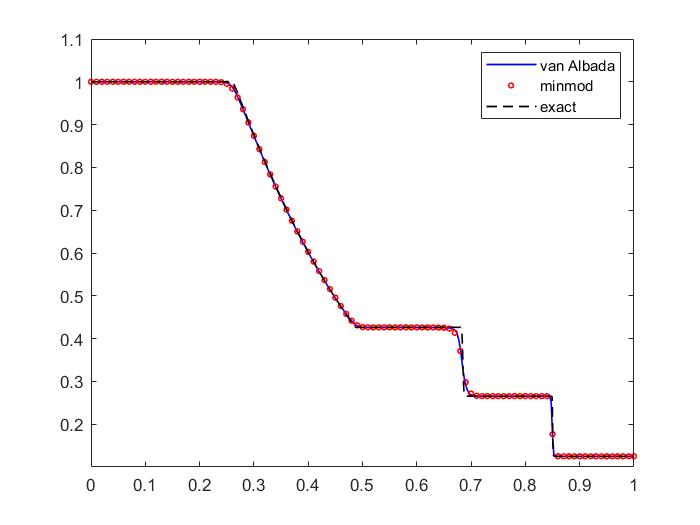}
    \caption{Sod's problem --- density computed with NT scheme and vA limiter, $N = 400$, $t = 0.2$.}
    \label{fig:shocktube NT}
\end{figure}

\begin{figure}[h!]
    \centering
    \includegraphics[scale = 0.4]{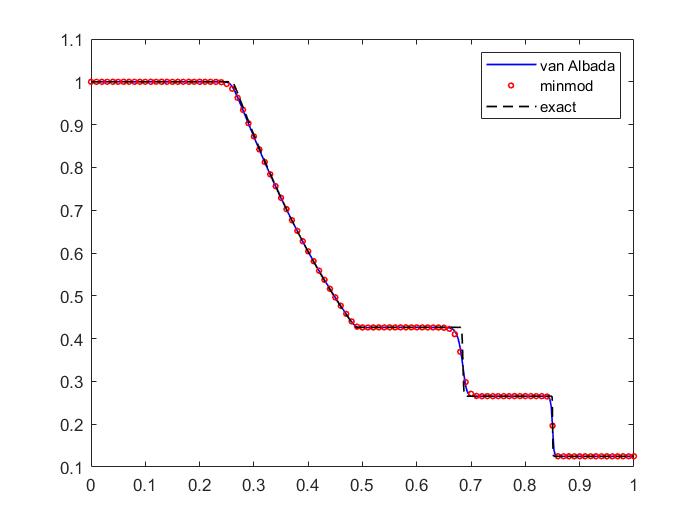}
    \caption{Sod's problem --- density computed with semi-discrete central-upwind scheme and vA limiter, $N = 400$, $t = 0.2$.}
    \label{fig:shocktube csup}
\end{figure}

As first test case,  we consider the Sod's shock tube problem: the computational domain is $[0, 1]$ with an interface at $x = 0.5$, and subject to initial states on the two sides of the interface  given by
\begin{equation*}
\left\{\begin{array}{c}
\begin{split}
        & p_{L} = 1.0, \quad p_{R} = 0.1  \\
        & \rho_{L} = 1.0, \quad \rho_{R} = 0.125 \\
        & u_{L} = 0.0, \quad u_{R} = 0.0 
\end{split}
\end{array}\right..
\end{equation*}
Figures \ref{fig:shocktube NT} and \ref{fig:shocktube csup} show the computed density at $t = 0.2$ with $N = 400$ mesh grids, computed with the NT and, respectively, semi-discrete central-upwind schemes. The numerical results, comparing the vA and minmod${}_1$ limiters show that the high resolution of the computation based on vA limiter  are on par with those based on the minmod limiter. In fact, we observe that the contact wave in Figure \ref{fig:shocktube NT} is even slightly better resolved by the NT+vA limiter compared with the central-upwind KNP+vA limiter in Figure \ref{fig:shocktube csup}. This agrees with the improved resolution of NT+vA limiter of linear advection discussed in remark \ref{rem:NTvA} and similar remark regarding the diffusion of the semi-discrete central-upwind KT schemes which goes back to \cite[p. 268]{KT00}.
\newline


\begin{figure}[h!]
    \centering
    \includegraphics[scale = 0.4]{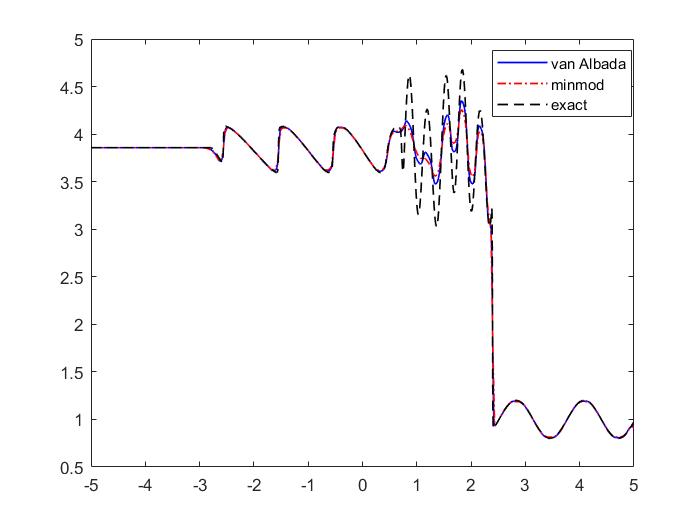}
    \caption{Osher-Shu problem --- density, semi-discrete central-upwind scheme, $N = 600$, $t = 1.8$.}
    \label{fig:shockentro csup}
\end{figure}

Next, we turn to a second test case of the Osher-Shu problem. The computational domain is $[-5 , 5]$ with the interface placed at $x = -4$, subject to initial data 
\begin{equation*}
\left\{\begin{array}{c}
\begin{split}
        & p_{L} = 10.33333, \quad p_{R} = 1.0,  \\
        & \rho_{L} = 3.857143, \quad \rho_{R} = 1+0.2\sin(5x), \\
        & u_{L} = 2.6293690, \quad u_{R} = 0.0. 
\end{split}
\end{array}\right.
\end{equation*}
In this problem, the solution of density consists of a discontinuity and a smooth harmonic wave. Figures \ref{fig:shocktube NT} and \ref{fig:shockentro csup} present the solutions of density given by the NT scheme and the semi-discrete central-upwind scheme computed with $N = 600$ mesh grids at $t = 1.8$. We observe  that the solutions computed with the vA limiter resolve the smooth extremum better than those computed with the minmod limiter; in particular, the solutions do not introduce spurious numerical oscillations at the location of the shock.

\begin{figure}[h!]
    \centering
    \includegraphics[scale = 0.4]{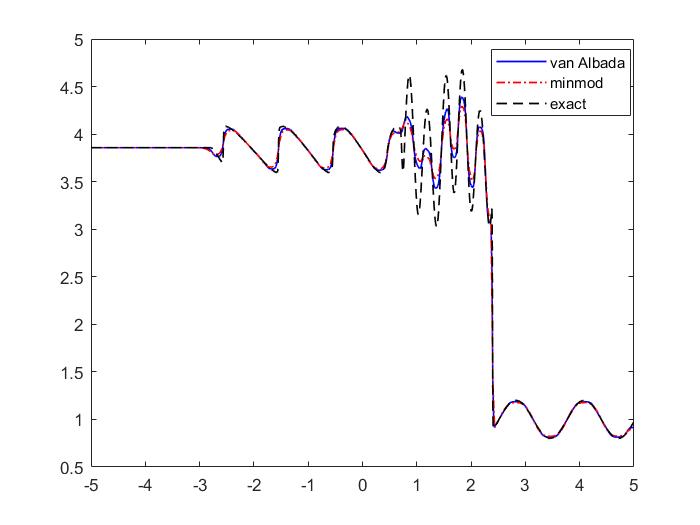}
    \caption{Osher-Shu problem --- density, NT scheme, $N = 600$, $t = 1.8$.}
    \label{fig:shockentro NT}
\end{figure}

\newpage

\subsection{Two-dimensional Euler equations}\label{sec:2DEuler}
\begin{figure}[h!]
    \centering
    \includegraphics[scale = 0.8]{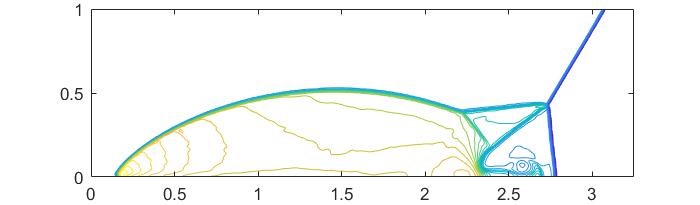}
    \caption{Double Mach reflection-density, 5th order WENO, $N_x = 480$, $N_y = 120$, $t = 0.2$.}
    \label{fig:DM WENO5}
\end{figure}

We examine the performance of the vA limiter in two-dimensional problems. We consider the double Mach reflection problem \cite{WC84} for the two-dimensional Euler equations,
\begin{equation*}
\frac{\partial}{\partial t}\left[\begin{array}{c}
     \rho  \\
     \rho u\\
     \rho v\\
     E
\end{array}\right]+
\frac{\partial}{\partial x}\left[\begin{array}{c}
     \rho u  \\
     \rho u^2+p\\
     \rho uv\\
     (E+p)u
\end{array}\right]+
\frac{\partial}{\partial y}\left[\begin{array}{c}
     \rho v  \\
     \rho uv\\
     \rho v^2+p\\
     (E+p)v
\end{array}\right] = 0, \quad p = (\gamma-1)\cdot[E-\frac{\rho}{2}(u^2+v^2)].
\end{equation*}
The problem is initiated with a Mach 10 oblique shock positioned at $(1/6,0)$ which makes a $60^{\circ}$ angel with $x-$axis. The computational domain is set to $[0,4]\times[0,1]$. The bottom boundary ($y=0$) consists of a reflecting wall beginning at $x = 1/6$, the short region from $x = 0$ to $x = 1/6$ is imposed with the initial post-shock conditions. The left boundary $(x = 0)$ is also assigned with the initial post-shock values, and at the right boundary $(x=4)$ all the gradients are set to zero. The values along the
top boundary $(y = 1)$ are set to describe the exact motion of the Mach 10 shock. We refer to \cite{WC84} for more detailed descriptions. We compute the solution of density at $t = 0.2$ with $480\times 120$ mesh grids. Figure \ref{fig:DM WENO5} shows the reference solution given by the 5th order WENO finite difference method \cite{Shu98}.

\begin{figure}[h!]
    \centering
    \includegraphics[scale = 0.8]{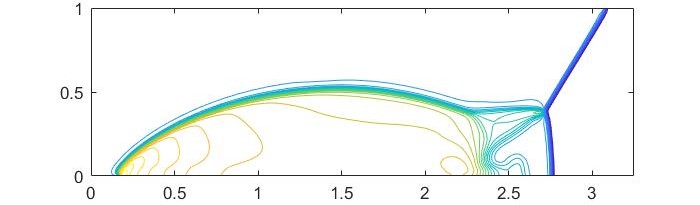}
    \caption{Double Mach reflection-density, NT scheme, minmod limiter, $N_x = 480$, $N_y = 120$, $t = 0.2$.}
    \label{fig:DM NT minmod}
\end{figure}

\begin{figure}[h!]
    \centering
    \includegraphics[scale = 0.8]{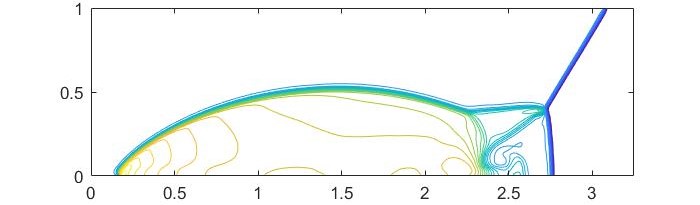}
    \caption{Double Mach reflection-density, NT scheme, vA limiter, $N_x = 480$, $N_y = 120$, $t = 0.2$.}
    \label{fig:DM NT VA}
\end{figure}

\begin{figure}[h!]
    \centering
    \includegraphics[scale = 0.8]{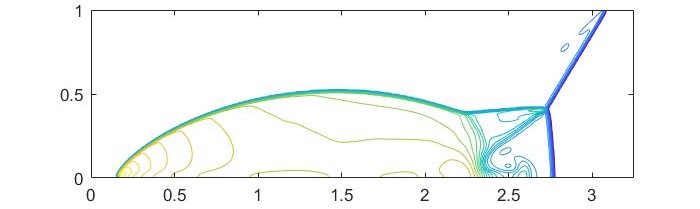}
    \caption{Double Mach reflection-density, semi-discrete central-upwind scheme, minmod limiter, $N_x = 480$, $N_y = 120$, $t = 0.2$.}
    \label{fig:DM csup minmod}
\end{figure}

\begin{figure}[h!]
    \centering
    \includegraphics[scale = 0.8]{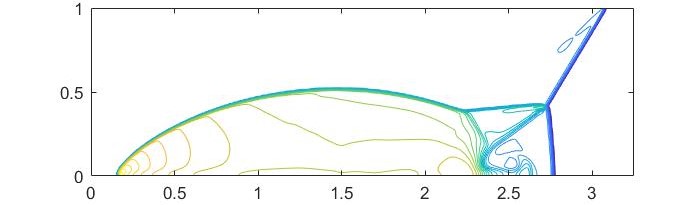}
    \caption{Double Mach reflection-density, semi-discrete central-upwind scheme, van Albada limiter, $N_x = 480$, $N_y = 120$, $t = 0.2$.}
    \label{fig:DM csup VA}
\end{figure}

Figures \ref{fig:DM NT minmod}, \ref{fig:DM NT VA} display the solutions computed with the NT scheme. It is clear from the results that the solution generated by the vA limiter has obviously better quality. The shock waves and the slip lines underneath are characterized by much sharper transitions. The computational results of the semi-discrete scheme are presented in Figures \ref{fig:DM csup minmod}, \ref{fig:DM csup VA}. Still, the numerical scheme captures more information of the near-wall flow when used in conjunction with the vA limiter. These results confirm the advantage of the vA limiter \eqref{eq:VA eps limiter} in terms of resolving multidimensional complex flow structures.

\vspace{0.5in}

\newpage

\subsection{Incompressible Euler equations}\label{sec:incompressible}
We consider the two-dimensional double shear layer problem for incompressible Euler equations,
\[
\left\{\begin{array}{l}
    \partial_t \bu+\nabla \cdot (\bu\otimes\bu) = -\nabla p,  \\ \\
     \nabla\cdot\bu = 0, 
\end{array}\right.
\]
where $p$ denotes the pressure and $\bu = (u,v)$ is the velocity field. The initial data is given by
\[
u(x,y,0) = \left\{\begin{array}{ll}
    \tanh{((y-\frac{\pi}{2})/\rho)} & y\leq\pi, \\ \\
    \tanh{((\frac{3\pi}{2}-y)/\rho)}  & y>\pi,
\end{array}\right. \qquad
v(x,y,0) = \delta \sin{(x)}
\]
over the rectangular domain $[0,2\pi]\times [0,2\pi]$. We take $\displaystyle \rho = \frac{\pi}{15}$ and $\delta = 0.05$. The equations are solved with the finite difference projection method (see e.g. \cite{LT97,KT97,Zhang09}), which consists of a predictor step and a corrector step
\[
\begin{split}
&\frac{\bu^*-\bu^n}{\Dt}+\nabla\cdot(\bu^n\otimes\bu^n) = 0,\\
& \frac{\bu^{n+1}-\bu^*}{\Dt}+\nabla p^{n+1} = 0 \quad \text{with}\quad  \nabla\cdot\bu^{n+1} = 0.
\end{split}
\]
The divergence free constraint, $\nabla\cdot\bu^{n+1} = 0$, leads to a Poisson equation for pressure
\[
-\Delta p^{n+1} = -\frac{1}{\Dt}\nabla\cdot\bu^*.
\]
In the predictor step, we reconstruct the convection fluxes, $\bu^n\otimes \bu^n$, with the Rusanov-type flux \eqref{eq:Rus flux} coupled with the second order MUSCL reconstructions \eqref{eq:MUSCL reconstruction}. In the corrector step, assuming the uniform mesh size is employed, $\Dx = \Dy = h$, the gradient, divergence and Laplace operators are approximated with the second order central difference operators,
\[
 \nabla^h \phi_{i,j}= \left(\begin{array}{c}
     \displaystyle \frac{\phi_{i+1,j}-\phi_{i-1,j}}{2\Dx}  \\
     \displaystyle \frac{\phi_{i,j+1}-\phi_{i,j+1}}{2\Dy}
\end{array}\right),
\]
\[
\nabla^h\cdot\bu_{i,j} =  \frac{u_{i+1,j}-u_{i-1,j}}{2\Dx}+\frac{v_{i+1,j}-v_{i-1,j}}{2\Dy}, \quad \bu = (u,v),
\]
\[
\Delta^h \phi_{i,j} = \frac{\phi_{i+2,j}-2\phi_{i,j}+\phi_{i-2,j}}{4\Dx^2}+ \frac{\phi_{i,j+2}-2\phi_{i,j}+\phi_{i,j-2}}{4\Dy^2}.  
\]
To achieve a higher order of accuracy in time, we apply the explicit SSPRK3 method for time marching. The time step is taken according to the CFL condition
\[
\Dt\max_{i,j}\Big(\frac{|u_{i,j}|}{\Dx}+\frac{|v_{i,j}|}{\Dy}\Big)\leq 0.5.
\]

\begin{figure}[h!]
    \centering
    \begin{subfigure}{0.4\textwidth}
     \includegraphics[scale = 0.3]{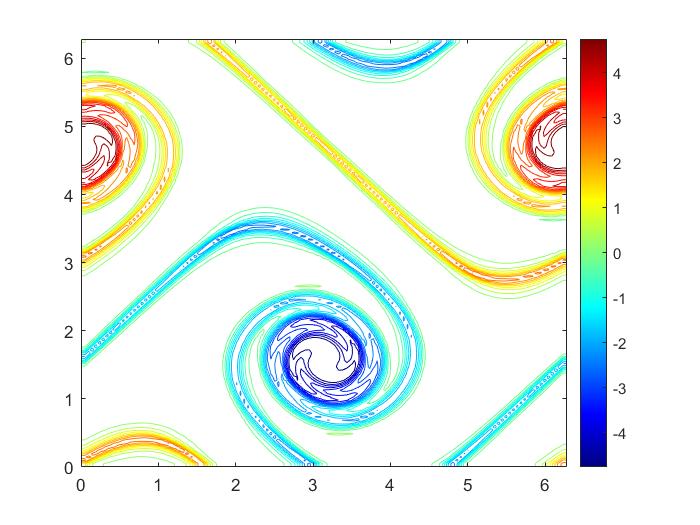}
     \subcaption{5th order WENO}
    \end{subfigure}
    \quad
    \begin{subfigure}{0.4\textwidth}
     \includegraphics[scale = 0.3]{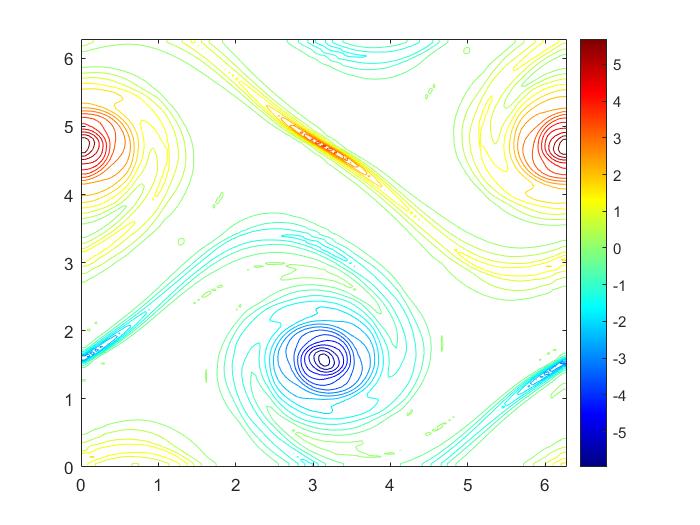}
     \subcaption{MUSCL-minmod}
    \end{subfigure}
    \\
    \begin{subfigure}{0.4\textwidth}
     \includegraphics[scale = 0.3]{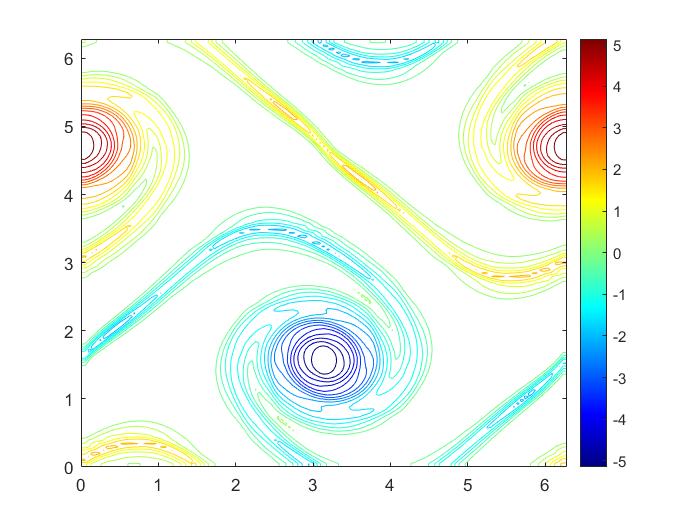}
     \subcaption{MUSCL-van Albada}
    \end{subfigure}
    \caption{Double shear layer, contours of vorticity. $N_x = N_y = 128$, $CFL = 0.5$, $t = 8$.}
    \label{fig:shear layer contours}
\end{figure}

\begin{figure}[h!]
    \centering
    \includegraphics[scale = 0.35]{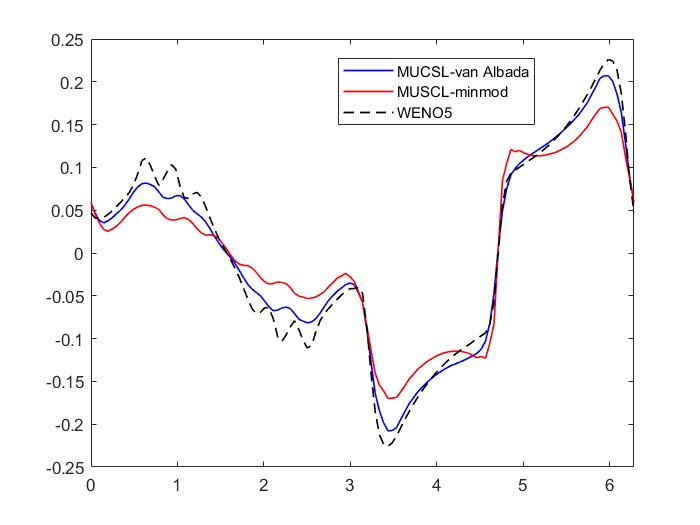}
    \caption{Double shear layer, solutions of $v$ along $x = \pi$, $N_x = N_y = 128$, $\CFL = 0.5$, $t = 8$.}
    \label{fig:shear layer velocity}
\end{figure}

Figure \ref{fig:shear layer contours} presents the solutions of vorticity at $t = 8$. We observe that the solution contours given by the vA limiter are more concentrated than those by the $\minmod$ limiter. Figure \ref{fig:shear layer velocity} compares the solutions of $v$ along the line $x=\pi$. The vA limiter shows clear advantage over the $\minmod$ limiter in terms of capturing sharp extrema.

\subsection{Concluding remarks}
The vA limiter dates back to the 1982 work \cite{ALR82}. Its TVD stability and related non-oscillatory properties were overshadowed by extensive studies of the class of minmod limiters, \cite{vLeer79,Harten83, vLeer84, Sweby84,CO85,OT88,NT90}. In this work, we revisit the vA limiter in the context of high-resolution central schemes, showing that its non-oscillatory properties   improve   corresponding computations using  the minmod limiter, or  are at least  on par with the minmod-based results in other computations.
This was demonstrated in the 1D Osher-Shu problem depicted in Figures \ref{fig:shockentro NT}--\ref{fig:shockentro csup} and more noticeably, in the 2D problems reported in sections \ref{sec:2DEuler} \& \ref{sec:incompressible}.
We do not claim that the use of vA limiter always outperforms the class of minmod limiters: the vA limiter are found on par with the minmod results for 1D Sod problem reported in Figures \ref{fig:shocktube NT}--\ref{fig:shocktube csup}; indeed, this is further highlighted in Figure \ref{fig:shocktube csupmm=1.2} below, where we zoom on the density variation of the Sod problem and observe comparable results for  the vA limiter which is compared with the less diffusive minmod${}_\theta$ with $\theta=\hf(1+\sqrt{2})$ (corresponding to the upper bound on the right of \eqref{eq:VA bound}). The main point here is to shed light on the smoother vA limiter, showing that it produces equal or even better results than the class of minmod limiters. The vA limiter was  added to the minmod limiter  in the Python code \cite{Python} (solver scheme `sd3'). 

\begin{figure}[h!]
    \centering
    \includegraphics[scale = 0.5]{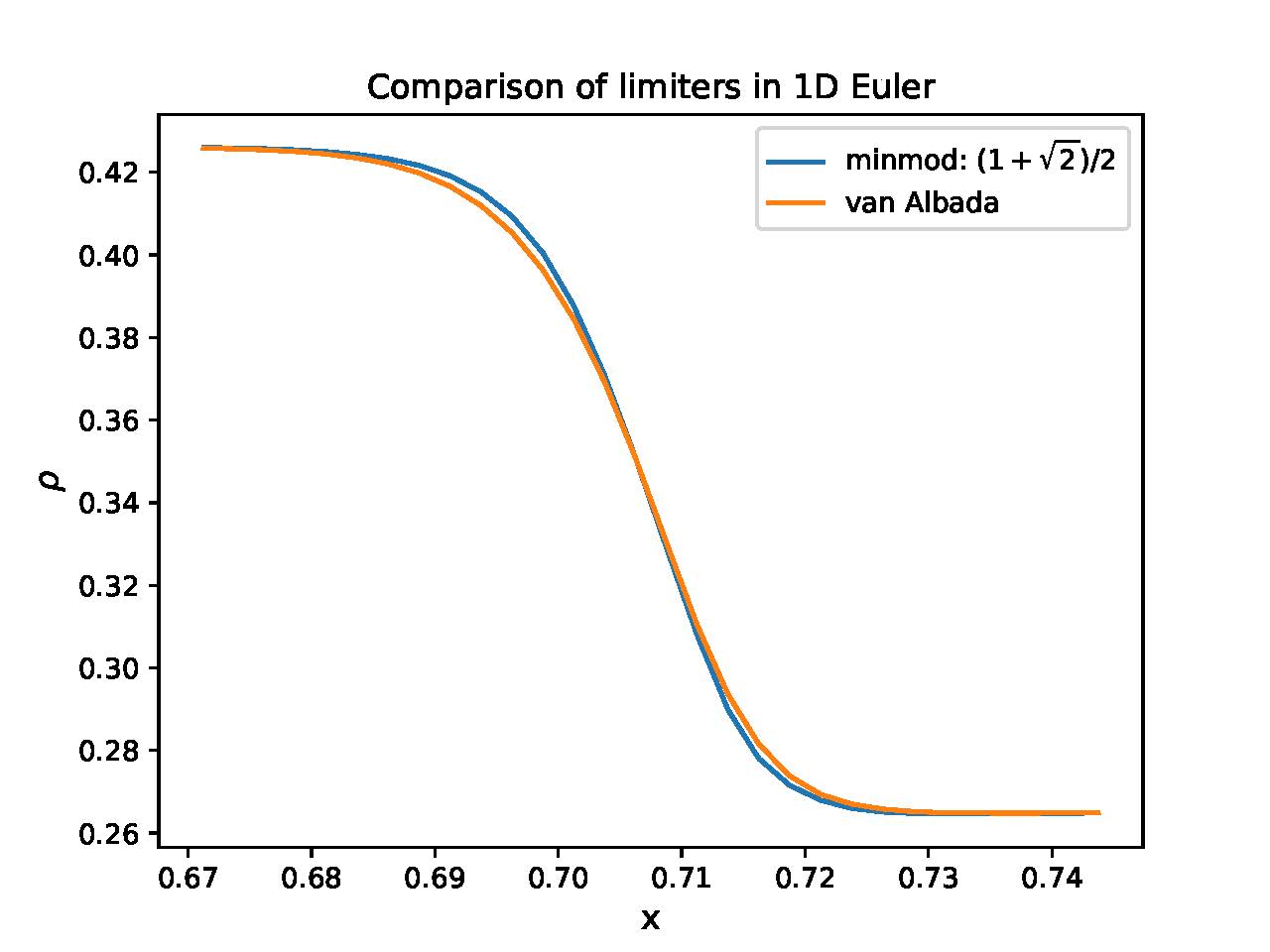}
    \caption{Sod's problem. Zoom on density computed by NT scheme, $N = 400$ at $t = 0.2$ using vA limter  vs. minmod${}_\theta$ limiter with $\theta=1.2$.}
    \label{fig:shocktube csupmm=1.2}
\end{figure}
\newpage

\bigskip\noindent
{\bf Statement}. On behalf of all authors, the corresponding author states  that this work is in compliance with Ethical Standards, that we did not participate in any Research involving Humans/and or Animals under Compliance with Ethical standards, and  that there is no conflict of interest.

\bibliographystyle{plain}
\bibliography{Ref}

\end{document}